\documentclass{article}
\usepackage{amsmath,amssymb,amsfonts,amsthm}
\makeindex

\title{Every set of first-order formulas is equivalent to an independent set\footnote{The following paper is a translation in English from the original paper of Iegor Reznikoff (\cite{Reznikoff}) ``Tout ensemble de formules de la logique classique est equivalent $\acute{a}$ un ensemble independant". This paper follows closely the arguments of Reznikoff in \cite{Reznikoff} , but it is not a word-by-word translation. It is intended only as a reference, not for publication. It is posted on arXiv with the permission of Dr. Reznikoff who we would like to thank. }}
\author{Ioannis Souldatos (for the translation)\\
Iegor Reznikoff (for the original French)}

\begin{document}
\maketitle

\newcommand{\omegaone}{\ensuremath{\omega_1}}
\newcommand{\lomegaone}{\ensuremath{\mathcal{L}_{\omega_1,\omega}}}
\newcommand{\alephs}[1]{\ensuremath{\aleph_{#1}}}
\newcommand{\alephalpha}{\alephs{\alpha}}
\newcommand{\alephomegaone}{\alephs{\omegaone}}
\newcommand{\alephalphaplus}{\alephs{\alpha+1}}
\newcommand{\beths}[1]{\ensuremath{\beth_{#1}}}
\newcommand{\bethalpha}{\beths{\alpha}}
\newcommand{\bethomegaone}{\beths{\omegaone}}
\newcommand{\bethalphaplus}{\beths{\alpha+1}}
\newcommand{\z}{\ensuremath{\mathcal{Z}}}
\newcommand{\M}{\ensuremath{\mathcal{M}}}
\newcommand{\N}{\ensuremath{\mathcal{N}}}
\newcommand{\A}{\ensuremath{\mathcal{A}}}
\newcommand{\B}{\ensuremath{\mathcal{B}}}
\newcommand{\F}{\ensuremath{\mathcal{F}}}
\newcommand{\D}{\ensuremath{\mathcal{D}}}
\newcommand{\C}{\ensuremath{\mathcal{C}}}
\newcommand{\E}{\ensuremath{\mathcal{E}}}
\newcommand{\G}{\ensuremath{\mathcal{G}}}
\newcommand{\mnneat}{(\M,\N)-neat\;}
\newcommand{\mnrich}{(\M,\N)-rich\;}
\newcommand{\mnfull}{(\M,\N)-full\;}
\newcommand{\mnkappafull}{$(\M,\N,\kappa)$-full\;}
\newcommand{\mneat}{\M-neat\;}
\newcommand{\mrich}{\M-rich\;}
\newcommand{\mfull}{\M-full\;}
\newcommand{\AminusMN}{\ensuremath{\A\setminus(\M\cup\N)\;}}
\newcommand{\lang}[1]{\ensuremath{\mathcal{L}_{#1}}}
\newcommand{\langhat}{\ensuremath{\widehat{\lang{}}}}
\newcommand{\phialpha}{\ensuremath{\phi_\alpha}}
\newcommand{\phikappa}{\ensuremath{\phi_\kappa}}
\newcommand{\phiM}{\ensuremath{\phi_{\M}}}
\newcommand{\philtok}{\ensuremath{\phi_{\ltok}}}
\newcommand{\phialphaplus}{\ensuremath{\phi_{\alpha+1}}}
\newcommand{\existslambda}{\ensuremath{\exists^{\omega\le\cdot\le\lambda}}}
\newcommand{\Malpha}{\ensuremath{\mathcal{M}_\alpha}}
\newcommand{\ch}{\ensuremath{\mathcal{CH}_{\omega_1,\omega}}}
\newcommand{\homch}{\ensuremath{\mathcal{HCH}_{\omega_1,\omega}}}
\newcommand{\ltok}{\ensuremath{\lambda^{\kappa}}}
\newcommand{\ktok}{\ensuremath{\kappa^{\kappa}}}
\newcommand{\ltoomega}{\ensuremath{\lambda^{\omega}}}
\newcommand{\lambdaalpha}{\ensuremath{\lambda_{\alpha}}}
\newcommand{\lambdaalphas}[1]{\ensuremath{\lambda_{\alpha_{#1}}}}
\newcommand{\dalpha}{\ensuremath{D_{\alpha}}}
\newcommand{\dalphaplus}{\ensuremath{D_{\alpha+1}}}
\newcommand{\dalphas}[1]{\ensuremath{D_{#1}}}
\newcommand{\taukappa}{\ensuremath{\tau_{\kappa}}}
\newcommand{\fuptox}[2]{\ensuremath{#1\upharpoonright_{(-\infty,#2)}}} 
\newcommand{\Sinfty}{\ensuremath{S_{\infty}}}
\newcommand{\Skappa}{\ensuremath{S_{\kappa}}}
\newcommand{\f}[1]{\ensuremath{f(\{#1\})}}
\newcommand{\Talpha}{T_\alpha}
\newcommand{\Dalpha}{D_\alpha}
\newcommand{\Tbeta}{T_\beta}
\newcommand{\phialphai}{\phi_{\alpha_i}}
\newcommand{\phibeta}{\phi_\beta}
\newcommand{\phigamma}{\phi_\gamma}
\newcommand{\psialpha}{\psi_\alpha}
\newcommand{\psialphai}{\psi_{\alpha_i}}
\newcommand{\psibeta}{\psi_\beta}
\newcommand{\psigamma}{\psi_\gamma}
\newtheorem{thrm}{Theorem}
\newtheorem{lem}[thrm]{Lemma}
\newtheorem{cor}[thrm]{Corollary}
\newtheorem{df}[thrm]{Definition}
\newtheorem{claim}{Claim}
\newtheorem{note}[thrm]{Note}

\abstract A set of first-order formulas, whatever the cardinality of
the set of symbols, is equivalent to an independent set.

In the following we work with classical (first-order) logic. The
Axiom of Choice is assumed.

\begin{df} Two sets of formulas are \emph{equivalent}, if any formula of
the one set is a consequence of the other and conversely. (Equiv.
they have the same models).

A set of formulas $T$ is \emph{independent}, if for all $\phi\in T$,
\[T\setminus\{\phi\}\nvDash\phi.\] (Equiv. there is a model for
$(T\setminus\{\phi\})\cup\{\neg\phi\}$).
\end{df}

\begin{thrm}\label{TarskisThrm}(Tarski) Every countable set of formulas is equivalent
to an independent set.
\begin{proof} Let $T=\{\phi_0,\phi_1,\ldots\}$ a countable set of
formulas. Without loss of generality there are no valid formulas in
$T$.

Define inductively
\begin{itemize}
  \item $\psi'_0=\phi_0$ and
  \item $\psi'_{n+1}=$ least $\phi_m$ such that
  $\psi'_0,\ldots,\psi'_n\nvDash \phi_m$.
\end{itemize}
It is not hard to see that $T$ is equivalent to the set
$\{\psi'_n|n\in\omega\}$. If this set is finite, then $T$ is
equivalent to its conjunction. So, assume it is infinite and define
\begin{itemize}
  \item $\psi_0=\psi'_0$ and
  \item $\psi_{n+1}=\bigwedge_{m\le n} \psi'_m\rightarrow
  \psi'_{n+1}$.
  \end{itemize}
Since $\psi'_0,\ldots,\psi'_n\nvDash \psi'_{n+1}$, there is a model
$\M$ that satisfies $\psi'_0,\ldots,\psi'_n$ and $\neg\psi'_{n+1}$.
Then $\M\nvDash \psi_{n+1}$, while $\M\models\psi_m$ for $m<n+1$.
For $m>n+1$, since $\M$ doesn't satisfy the antecedent of $\psi_m$,
it trivially satisfies $\psi_m$. Therefore
\[\M\models\bigwedge_{m\neq n+1}\psi_m\wedge\neg\psi_{n+1},\]
witnessing the fact that $\{\psi_n|n\in\omega\}$ is an independent
set.

Moreover, it is an easy induction to see that the sets
$\{\psi'_n|n\in\omega\}$ and $\{\psi_n|n\in\omega\}$ are equivalent,
which finishes the proof.
\end{proof}
\end{thrm}

\begin{lem}\label{ReznikoffsLemma}Let $C,D$ be two disjoint sets such that:
\begin{itemize}
  \item $|D|\le|C|$ and
  \item For all $\phi\in C$, $(C\cup D)\setminus\{\phi\}\nvDash \phi$.
  \end{itemize}
  Then $C\cup D$ is equivalent to an independent set.
\begin{proof}Let $f$ be an injection from $D$ to $C$. Then
\[\{\psi\wedge f(\psi)|\psi\in D\}\cup (C\setminus f(D))\] is an independent
set equivalent to $C\cup D$.
\end{proof}
\end{lem}

Now, let $T$ be a set of formulas and without loss of generality,
there are no valid formulas in $T$ (valid formulas are equivalent to
the empty set). For a formula $\phi\in T$, denote by $S(\phi)$ the
set of symbols that appear in $\phi$ and let
\[S=\bigcup_{\phi\in T} S(\phi).\]
Without loss of generality $S$ is infinite. Otherwise $T$ would be at
most countable and equivalent to an independent set by Theorem
\ref{TarskisThrm}. If $S$ is infinite, then $S$ and $T$ have the same
cardinality and let \[|S|=|T|=\kappa\ge\omega.\]

We partition $T$ into sets $\Talpha$, $\alpha<\kappa$ as follows:

For $\alpha=0$, fix a formula $\phi_0\in T$ and let $T_0=\{\psi\in
T| S(\psi)\subset S(\phi_0)\}$. For $0<\alpha<\kappa$, assume that
we have defined $\phibeta$ and $\Tbeta$, for all $\beta<\alpha$. By
a cardinality argument, \[S\setminus\bigcup_{\beta<\alpha}
S(\phibeta)\neq\emptyset.\] Therefore, there exists a formula
$\phialpha$ that contains a symbols that doesn't appear in any of
the $\phibeta$, $\beta<\alpha$. Define
\[N_\alpha=S(\phialpha)\setminus\bigcup_{\beta<\alpha}
S(\phibeta),\]the set of new symbols that appear in $\phialpha$.
Then $N_\alpha\neq\emptyset$ and define
\[\Talpha=\{\psi\in T|S(\psi)\subset\bigcup_{\beta\le\alpha} S(\phibeta)\mbox{ and } S(\psi)\cap N_\alpha\neq\emptyset\},\]
i.e. $\Talpha$ is the set of formulas in which appears one of the
new symbols in $N_\alpha$.

Then $T=\bigcup_{\alpha<\kappa} \Talpha$ and the different
$\Talpha$'s are disjoint.

\begin{df} If $\psi\in \Talpha$ and $S(\psi)\cap N_\beta\neq\emptyset$, for
$\beta\le\alpha$, denote this by $\beta|\psi$. In particular, for
$\psi\in \Talpha$, $\alpha|\psi$.

If $\beta|\phialpha$, with $\beta<\alpha$, denote this by
$\beta||\phialpha$.
\end{df}
Observe here that the first definition is for any $\psi\in T$, while
the second one is only for the $\phialpha$'s. Also, for any $\psi\in
T$, there are only finitely many $\beta$'s with $\beta|\psi$.

Now let \[\psialpha=\bigwedge_{\beta||\phialpha} \phibeta\rightarrow
\phialpha,\] if there exists such a $\beta$. Otherwise, let
$\psialpha=\phialpha$. Denote by $C$ the set of all the
$\psialpha$'s.

On the other hand, for $\phi\neq\phialpha$, all $\alpha<\kappa$, let
\[\phi'=\bigwedge_{\beta|\phi}\phibeta\rightarrow\phi.\] As we noted, there is always such a $\beta$. Denote
\[\Dalpha=\{\phi'=\bigwedge_{\beta|\phi}\phibeta\rightarrow\phi|\phi\in \Talpha\mbox{ and } \phi\neq\phialpha\}\]
and let $D=\bigcup_{\alpha}\Dalpha$. ($D$ may be empty. We can not
exclude this possibility).

\begin{lem} Suppose that $T$ satisfies the following condition:
\[(\star)\hspace{.5cm}\mbox{If } \psi,\phi_1,\ldots,\phi_n\in T \mbox{ and }S(\psi)\nsubseteq\bigcup_{i=1}^n S(\phi_i), \mbox{ then }
\{\phi_1,\ldots,\phi_n\}\nvDash \psi.\] Then $C$ and $D$ as defined
above, satisfy the conditions of Lemma \ref{ReznikoffsLemma} and $T$
is equivalent to an independent set.
\begin{proof}First of all it is clear that $|C|=\kappa\ge|D|$. It
also follows easily by induction on $\alpha<\kappa$ that the set
$\bigcup_{\beta<\alpha} \Tbeta$ is equivalent to the set
$\bigcup_{\beta<\alpha} (\{\psibeta\}\cup D_\beta)$. This implies
that $T$ is equivalent to $C\cup D$ and it suffices to verify that
for $\psialpha\in C$, $\psialpha$ is not a consequence of the other
elements of $C\cup D$:

Let $\psialpha=\bigwedge_{\beta||\phialpha} \phibeta\rightarrow
\phialpha$. Then the elements of $C\cup D$ different than
$\psialpha$ are of the form $\psigamma=\bigwedge_{\beta||\phigamma}
\phibeta\rightarrow \phigamma$, with $\gamma\neq\alpha$, or of the
form $\phi'=\bigwedge_{\beta|\phi}\phibeta\rightarrow\phi$, with
$\phi\neq\psialpha$, for all $\alpha<\kappa$.

Consider the implication
\[\models\left(\bigwedge_{i=1 \\ \alpha_i\neq\alpha}^m\psialphai\bigwedge_{j=1}^n \phi'_j\right)\rightarrow\psialpha.\]

Assume that $\psi_{\alpha_1},\ldots,\psi_{\alpha_m}\neq\psialpha$ and
$\alpha\nmid\phialphai$, for $i=1,\ldots,p$, while
$\alpha|\phialphai$, for $i=p+1,\ldots,m$. Similarly, assume that
$\phi'_1,\ldots,\phi'_q$ are such that $\alpha\nmid\phi_j$,
$j=1,\ldots,q$, while for $\phi'_{q+1},\ldots,\phi'_n$,
$\alpha|\phi_j$, $j=q+1,\ldots,n$.

Then \[S(\phialpha)\nsubseteq\bigcup_{i=1}^p S(\phialphai) \mbox{
and  } S(\phialpha)\nsubseteq\bigcup_{j=1}^q S(\phi'_j).\] Also, by
the definition of $\phialpha$,
\[ S(\phialpha)\nsubseteq\bigcup_{\beta||\phialpha} S(\phibeta).\]

By $(\star)$, there is a model $\M$ in which $\phialpha$ is false,
while all of the $\phi_{\alpha_1},\ldots,\phi_{\alpha_p}$,
$\phi_1,\ldots,\phi_q$ and $\{\phibeta|\;\beta||\phialpha\}$ are
true. Then $\psialpha$ is false in $\M$, while
$\psi_{\alpha_1},\ldots,\psi_{\alpha_p},\phi'_1,\ldots,\phi'_q$ are
true.

In addition, for $i=p+1,\ldots,m$, $\phialpha$ is among the
$\phigamma$'s in the conjunction of
$\psialphai=\bigwedge_{\gamma||\phialphai} \phigamma\rightarrow
\phialphai$ and the same is true for the conjunction of
$\phi'_j=\bigwedge_{\gamma|\phi_j} \phigamma\rightarrow \phi_j$, for
$j=q+1,\ldots,n$. Since $\phialpha$ is false in $\M$, then
$\psi_{\alpha_{p+1}},\ldots,\psi_{\alpha_m},\phi'_{q+1},\ldots,\phi'_n$
are trivially true. This proves that there is a model $\M$ that satisfies
all the $\psi_{\alpha_1},\ldots,\psi_{\alpha_m},\phi'_1,\ldots,\phi'_n$,
but which does not satisfy $\psialpha$. In other words, $\psialpha$ can
not be a consequence of other elements of $C\cup D$.
\end{proof}
\end{lem}
What remains is to prove that $T$ can be taken to satisfy $(\star)$. We use
Craig's Interpolation Theorem which we mention without proof.
\begin{thrm}(Craig) If $\psi\models\phi$, then there is a formula
$\tau$ such that
\begin{itemize}
  \item $\psi\models\tau$ and $\tau\models\phi$, and
  \item the non-logical symbols of $\tau$ appears in both $\psi$ and
  $\phi$.
\end{itemize}
$\tau$ is called the \emph{interpolant} between $\psi$ and $\phi$.
\end{thrm}

\begin{lem} Every set of non-valid formulas $T$ is equivalent to
a set of formulas that satisfies $(\star)$.
\begin{proof} Let \[E_1=\{\phi|\;T\models\phi\mbox{ and }
|S(\phi)|=1\}\] and
\[E_n=\{\phi|\;T\models\phi,\bigcup_{m<n} E_m\nvDash\phi\mbox{ and }
|S(\phi)|=n\}.\]

It is immediate that $T'=\cup_n E_n$ is equivalent to $T$. Let
$\psi,\phi_1,\ldots,\phi_n\in T'$ such that
$S(\psi)\nsubseteq\bigcup_{i=1}^n S(\phi_i)$. If we assume that
\[\{\phi_1,\ldots,\phi_n\}\models\psi,\] then by Craig's
Interpolation Theorem, there is a $\tau$ such that
\begin{itemize}
  \item $\{\phi_1,\ldots,\phi_n\}\models\tau$ and $\tau\models\psi$, and
  \item $S(\tau)\subset S(\psi)\cap (\cup_{i=1}^n S(\phi_i))$.
\end{itemize}
By the assumption on $\psi$, it must be $S(\tau)\subsetneq S(\psi)$
and $\psi\in T'$ would be a consequence of $\tau$ with $T\models
\tau$ and $|S(\tau)|<|S(\psi)|$, contradicting the definition of
$T'$.

Therefore, $T'$ satisfies $(\star)$.
\end{proof}
\end{lem}

Putting all the previous lemmas together we conclude
\begin{thrm}(Reznikoff) Every set of formulas is equivalent to an
independent set.
\end{thrm}

\end{document}